\documentclass[10pt]{article}

\usepackage{psfrag}
\usepackage{natbib}
\usepackage{enumerate}
\usepackage{array}
\usepackage{mathrsfs} 
\usepackage{amsfonts, amssymb, latexsym, slashbox}
\usepackage{amsmath,amsthm}
\usepackage{cancel}
\usepackage{url, graphicx, graphics, bm}
\usepackage{hyperref}
\usepackage{changes}  
\usepackage{comment}
\def\em{\it}
\def\emph{\textit}
\input{colordvi}

\def\o*{o_{{\!}_{P^*}}}
\def\O*{\cal O_{{\!}_{P^*}}}

\def\E{\mathrm{E}}
\def\Var{\mathrm{Var}}

\def\P{\mathrm{Pr}}

\def\P{\mathrm{P}}

\def\.{\mbox{.}}

\def\argmax{\mathop{\mbox{argmax}}}

\def\non{\nonumber}

\def\diag{\mbox{diag}}

\def\sfrac(#1,#2){\mbox{$\frac{#1}{#2}$}}
\def\ints(#1,#2){\mathbb{I}_{#1}^{#2}}

\def\T{{ \mathrm{\scriptscriptstyle T} }}

\newtheorem {theorem}{Theorem}

\newtheorem {rem}{Remark}
\newtheorem{assump}{}

\newtheorem{prop}[theorem]{Proposition}
\newtheorem{condition}{Condition}
\title{Maximum Approximate Likelihood Estimation in Accelerated Failure Time Model for Interval-Censored Data} 
\author{Zhong Guan\\
Department of Mathematical Sciences\\
Indiana University South Bend, USA\\
zguan@iusb.edu\\
}
\begin{document}
\maketitle

\section*{Abstract}
The approximate Bernstein polynomial model, a mixture of beta distributions, is applied to obtain maximum likelihood estimates of the regression coefficients, and the baseline density and survival functions in an accelerated failure time model based on interval censored data including current status data. The rate of convergence of the proposed estimates are given under some conditions for uncensored and interval censored data. Simulation shows that the proposed method is better than its competitors.
The proposed method is illustrated by fitting the Breast Cosmetic Data using the accelerated failure time model.

\noindent\textbf{Key Words and Phrases:} Accelerated Failure Time Model;  Approximate Likelihood;
Beta mixture model; Current Status Data; Interval Censoring; Smooth Estimation; Survival Curve.
\section{Introduction}
When a model in statistics involves some infinite-dimensional parameters such as a totally unspecified underlying distribution, it should rather be called a {\em probability model} or a non- or semi-parametric {\em statistical problem} to be not confused with a working statistical model which is used to solve the problem.  Because the sample size is always finite a working model must have a finite dimensionality which can be unknown. The traditional parametric models are of known dimensions while many mixture models are of unknown dimensions.

Due to the lack of an appropriate approximate model for the unspecified underlying baseline distribution, it is much more difficult to estimate the AFT model than PH and PO models using maximum likelihood method based on interval censored data. Traditionally we use step-function to approximate an unknown smooth distribution function so that we have a finite-dimensional working model which is (discrete) multinomial model and results in empirical distribution, Kaplan-Meier estimator \citep{Kaplan-Meier-1958-JASA}, Turnbull estimator \citep{Turnbull-1976-JRSSB}, and empirical likelihood method \citep{Owen1988}, and so on. This works because  the step-functions are dense in the space of certain continuous functions. Despite the roughness of resulting maximum likelihood estimate, this approximate model works quite well for many complicated statistical problems including the analysis of incompletely observed data such as data containing censored, grouped, truncated, and even missing values. An important example is the analysis of interval-censored event time data using PH \citep{Cox1972} and PO \citep{Pettitt-JRSSB-1982,Bennett-stats-Med-1983} models although the  semiparametric maximum likelihood estimate is not necessarily unique. However,  if data are interval-censored it is impossible to find a semiparametric maximum likelihood estimate of the baseline distribution using this approach for AFT model \citep{Kalbfleisch-Prentice-1980-book}, an important alternative to the PH  and PO models.
To the knowledge of the author, most inference procedures so far for AFT model with unspecified baseline distribution focus on the estimation of the regression coefficients \citep{Tian-and-Cai-bka-2006} and the right-censored \citep{Buckley-and-James-1979-bka,Tsiatis-1990-aos,Wei-etal-1990-bka,Jin-etal-2003-bka} or the current status data in biostatistics \citep{Huang-and-Wellner-1997} or the binary choice model in econometrics \citep{Cosslett-1983-Ecometrica,Cosslett-1987-Ecometrica,Cosslett-2004-Ecometrica, Klein-and-Spady-1993-econometric}.

A few most relevant works to this paper include
\cite{Hanson-and-Johnson-2004-jcgs-bayesian}, \cite{Komarek-etal-jcgs-2005-AFT-spilne}, and \cite{Zhang-and-Davidian-2008-biom-Smooth-SP}. The first proposed a Bayesian semiparametric AFT model for estimating survival and density functions, the second used B-splines with penalties to smooth the error density  with some candidate parametric models, and the latter proposed smooth estimates of survival function for PH, PO and AFT models using the so-called  seminonparametric (SNP) density \citep{Gallant-and-Nychka-1987-econom} which is a truncated Hermite series approximation of a density function. However Bernstein polynomials seem much better dense functions than step-functions and others for the purpose of building working statistical models \citep{Guan-2017-jns}. This Bernstein polynomial approximation is actually a mixture of some specified beta distributions with shapes related to the degree. This model has been successfully applied to grouped, contaminated, multivariate, and interval censored data \citep{Guan-2017-jns,Guan-2019-mable-deconvolution,Guan-arXiv-2019, Wang-and-Guan-2019}.
This model shall be applied to find maximum likelihood estimates of the regression coefficients, and the density and survival functions in the AFT model.


 \section{Methodology}
Let $T$ be an event time and $\bm X$ be an associated $d$-dimensional covariate with distribution $H(\bm x)$ on $\cal{X}$.
Let $f(t\mid\bm x)$ and $S(t\mid\bm x)$ be, respectively, the density and survival functions of failure time $T$ given $\bm X=\bm x$.
The AFT model can be specified as
\begin{equation}\label{eq: AFT model in terms of density}
f(t\mid\bm x)=f(t\mid\bm x;\bm\gamma)=e^{-\bm\gamma^\T\bm x}f(te^{-\bm\gamma^\T\bm x}\mid\bm 0),\quad t\in[0,\infty),
\end{equation}
where $\bm\gamma\in \mathbb{G}\subset \mathbb{R}^d$. Let $\bm\gamma_0\in \mathbb{G}$ be the true value of $\bm\gamma$.
The AFT model (\ref{eq: AFT model in terms of density}) is equivalent to
$$S(t\mid\bm x;\bm\gamma)=S(te^{-\bm\gamma^\T\bm x}\mid\bm 0),\quad t\in[0,\infty).$$
Thus this is actually a scale regression model.
The AFT model can also be written as linear regression $\log(T)=\bm\gamma^\T\bm x +\varepsilon$. It is clear that one can choose any $\bm x_0$ in $\mathcal{X}$ as baseline by transform $\tilde{\bm x}=\bm x-\bm x_0$. If $f(t\mid \bm 0)$ has support $[0,\tau_0)$, $\tau_0\le\infty$, then 
$f(t\mid \bm x)$ has support $[0,\tau_0 e^{\bm\gamma_0^\T\bm x})$.  We define $\tau=\max\{\tau_0 e^{\bm\gamma_0^\T\bm x}: \bm x\in\mathcal{X}\}$ if $\tau_0<\infty$ and $\tau=\infty$ otherwise.
The above AFT model can also be written as
$$f(t\mid\bm x;\bm\gamma)=e^{-\bm\gamma^\T {\bm x}} f_0(te^{-\bm\gamma^\T {\bm x}}), \quad
S(t\mid\bm x;\bm\gamma)=S_0(te^{-\bm\gamma^\T {\bm x}}),$$
where   $f_0(t)=f(t\mid\bm  0)$ and $S_0(t)=S(t\mid\bm 0)=\int_t^\infty f_0(u)du$.
Clearly, the above model is also true for the transformed failure time $T^*=aT$ for any  $a>0$.

With interval censoring, the observable random variables are $\bm Z=(\Delta, \bm X, \bm Y)$, where $\bm Y=(Y_1,Y_2]$ and $\Delta$ is the censoring indicator, i.e.,  uncensored $T=Y=Y_1=Y_2$ if $\Delta=0$, and  interval  censored $T\in\bm Y=(Y_1,Y_2]$, $0\le Y_1<Y_2\le\infty$, if $\Delta=1$. For Case 1 interval censored data, i.e., the current status or doubly censored data, $\bm Y=(0,U]$ or $\bm Y=(U,\infty)$. In this case the distribution function of the examination time $U$ given $\bm X=\bm x$ is denoted by $G(u\mid \bm x)$. As in \cite{Huang-and-Wellner-1997} we reduce the cases with more than two examination times to the case with two examination times, i.e., the Case $2$ interval censored data, and denote the joint distribution function of the observed examination times $\bm U=(U_1, U_2)$ given $\bm X=\bm x$  by $G(\bm u\mid \bm x)$.

For an observation $\bm z=(\delta, \bm x, \bm y=(y_{1},y_{2}])$, the exact full loglikelihood, up to an additive term independent of $(\bm\gamma, f_0)$, is
\begin{align}\nonumber
\ell(\bm\gamma, f_0;\bm z)&=(1-\delta)\log f(y\mid\bm x;\bm\gamma)
+\delta \log\{S(y_{1}\mid\bm x;\bm\gamma)-S(y_{2}\mid\bm x;\bm\gamma)\}\\\nonumber
&=(1-\delta)\{-\bm\gamma^\T {\bm x} +\log f_0(ye^{-\bm\gamma^\T {\bm x}})\}
\\\label{eq: likelihood for AFT model}
&~~~
+\delta \log\{S_0(y_{1}e^{-\bm\gamma^\T {\bm x}})-S_0(y_{2}e^{-\bm\gamma^\T {\bm x}})\}.
\end{align}
Let $\bm z_i=(\delta_i,\bm x_i, \bm y_i=(y_{1i},y_{2i}])$, $i=1,\ldots,n$, be independent observations of $\bm Z$.
The loglikelihood of the data is $\ell(\bm\gamma,f_0) =\sum_{i=1}^n \ell(\bm\gamma, f_0;\bm z_i)$. The Hessian matrix is $H_n(\bm\gamma,f_0)= {\partial^2\ell(\bm\gamma,f_0)}/{\partial \bm\gamma\partial \bm\gamma^\T}=\sum_{i=1}^n {\partial^2\ell(\bm\gamma,f_0;\bm z_i)}/{\partial \bm\gamma\partial \bm\gamma^\T}$.

The exact full likelihood $\ell(\bm\gamma,f_0)$ cannot be maximized unless $f_0$ and $S_0$ are specified. Because $y$ and $\bm x$ cannot be separated, step-function approximation does not work and it is also impossible to obtain a partial likelihood as that of \cite{Cox1972}.

In the case where $\tau_0=\infty$ or $\tau_0$ unknown   we choose $\tau_n>y_{(n)}= \max\{y_{i1}, y_{j2}: y_{j2}<\infty;\, i,j=1,\dots,n\}$ so that $S(\tau_n)$ and $\max_{\bm x\in\mathcal{X}} S(\tau_n\mid \bm x)$ are believed very small.
Then we approximate $f_0(t)$ and $S_0(t)$ on $[0,\tau_n]$, respectively, by
\begin{align*}
f_0(t)&\approx f_m(t; \bm p)=\frac{1}{\tau_n}\sum_{j=0}^m p_j\beta_{mj}\Big(\frac{t}{\tau_n}\Big),\quad t\in[0,\tau_n];\\
S_0(t)&\approx S_m(t; \bm p)= \sum_{j=0}^{m} p_j \bar B_{mj}\Big(\frac{t}{\tau_n}\Big),\quad  t\in[0,\tau_n],
\end{align*}
where
 $\beta_{mj}(t)=(m+1){m\choose  j}t^j(1-t)^{m-j}$,  $\bar B_{mj}(t)=1-\int_0^t\beta_{mj}(u)du$, $j=0,\dots,m$,   $S_m(\infty; \bm p)=0$, and $\bm p= (p_0,\ldots,p_{m})^\T$  satisfies constraint
\begin{equation}\label{eq: constr for p=(p0,...,p_{m})}
\bm p\in \mathbb{S}_{m}\equiv\Big\{\bm u=(u_0,\ldots,u_m)^\T: u_j\ge 0,  \; j=0,\dots,m, \;
 \sum_{j=0}^{m}u_j=1\Big\}.
\end{equation}
%
  Then $f(t\mid\bm x; \bm\gamma)$ and $S(t\mid\bm x; \bm\gamma)$  can be approximated, respectively, by
\begin{align}\non
f_m(t\mid\bm x; \bm\gamma, \bm p)&=e^{-\bm\gamma^\T {\bm x}}f_m\Big(te^{-\bm\gamma^\T {\bm x}};\bm p\Big)\\
&=
 \frac{e^{-\bm\gamma^\T {\bm x}}}{\tau_n}\sum_{j=0}^m p_j\beta_{mj}\Big(e^{-\bm\gamma^\T {\bm x}}\frac{t}{\tau_n}\Big),\quad t\in[0,\tau_ne^{\bm\gamma^\T {\bm x}}];\\\non
S_m(t\mid\bm x; \bm\gamma, \bm p)&=S_m\Big(te^{-\bm\gamma^\T {\bm x}};\bm p\Big)\\
&=  \sum_{j=0}^{m} p_j\bar {\mathcal{B}}_{mj}\Big(e^{-\bm\gamma^\T {\bm x}}\frac{t}{\tau_n}\Big),\quad t\in[0,\tau_ne^{\bm\gamma^\T {\bm x}}].
\end{align}
The likelihood $\ell(\bm\gamma,f_0)$ can be approximated by $\ell_m(\bm\gamma,\bm p)=\sum_{i=1}^n \ell_m(\bm\gamma,\bm p;\bm z_i),
$
where
\begin{align}\nonumber
\ell_m(\bm\gamma,\bm p;\bm z)
&= (1-\delta)\Big\{{-\bm\gamma^\T {\bm x}}+\log \sum_{j=0}^m p_j\beta_{mj} (ye^{-\bm\gamma^\T {\bm x}}/\tau_n )-\log\tau_n\Big\}\\\label{eq: approx bernstein loglikelihood for AFT model}
&~~~ +\delta \log\{S_m(y_{1}\mid\bm x; \bm\gamma, \bm p)-S_m(y_{2}\mid\bm x; \bm\gamma, \bm p)\}.
 \end{align}
 If $\tau_0$ is known, we choose $\tau_n=\tau_0$. If data are right-censored then $y_2=\tau_n$.
If $\tau_n\ne 1$ we divide all the observed times by $\tau_n$.  Thus we assume in the rest of this section that $\tau_n=1$.
%

 For a given degree $m$, let $(\hat{\bm\gamma}, \hat{\bm p})$ be a maximizer of $\ell_m(\bm\gamma,\bm p)$.
The change-point method \citep{Guan-jns-2015} applies for finding an optimal degree $m$.
%
For each $i=0,\dots,k$, fit the data to obtain $(\hat{\bm\gamma},\hat{\bm p})$ and $\ell_i=  \ell_{m_i}(\hat{\bm\gamma},\hat{\bm p})$, where $m_i=m_0+i$.
The optimal degree is $\hat m=\min\argmax_{1\le i\le k}\{R(m_i)\}$, where $R(m_i)=k\log\{({\ell_k-\ell_0})/{k}\}-i\log\{({\ell_i-\ell_0})/{i}\}-(k-i)\log\{({\ell_k-\ell_i})/({k-i})\}$, $i=1,\dots,k-1$ and $R(m_k)=0$.
With an optimal degree $m=\hat m$, $\hat{\bm\theta}=(\hat{\bm\gamma}, \hat{\bm p})$ is called a maximum approximate Bernstein likelihood estimator (MABLE) of $\bm\theta=(\bm\gamma,\bm p)$. The resulting MABLEs of $f(t\mid\bm x)$ and $S(t\mid\bm x)$ are, respectively,
 $\hat f_\mathrm{B}(t\mid\bm x)=f_m(t\mid \bm x; \hat{\bm\gamma}, \hat{\bm p})$ and $\hat S_\mathrm{B}(t\mid\bm x)=S_m(t\mid \bm x; \hat{\bm\gamma}, \hat{\bm p})$. The variance-covariance matrix of $\hat{\bm\gamma}$ can be estimated by $\hat\Sigma_{\bm\gamma}=-n\{H_n(\hat{\bm\gamma}; f_m(\cdot;\hat{\bm p}))\}^{-1}$.

 The derivatives of $\ell_m(\bm\gamma,\bm p;\bm z)$ with respect to $\bm p$ are
 \begin{align}\nonumber\label{eq: 1st deriv of ell_m wrt p for AFT model}
\frac{\partial\ell_m(\bm\gamma,\bm p;\bm z)}{\partial \bm p}
&=\frac{(1-\delta)\bm\beta_{m}(ye^{-\bm\gamma^\T {\bm x}})}{\sum_{j=0}^m p_j\beta_{mj}(ye^{-\bm\gamma^\T {\bm x}})}
 \\&~~~~~~+\frac{\delta \{\bar{\bm B}_m( y_1e^{-\bm\gamma^\T {\bm x}} )-\bar{\bm B}_m(y_2 e^{-\bm\gamma^\T {\bm x}} )\}}{S_m(y_{1}\mid\bm x; \bm\gamma, \bm p)-S_m(y_{2}\mid\bm x; \bm\gamma, \bm p)},
 \\\nonumber\label{eq: 2nd deriv of ell_m wrt p for AFT model}
\frac{\partial^2\ell_m(\bm\gamma,\bm p;\bm z)}{\partial \bm p\partial \bm p^\T}&=- \frac{(1-\delta) \big\{\bm\beta_{m}(ye^{-\bm\gamma^\T {\bm x}})\big\}^{\otimes 2}}{\{\sum_{j=0}^m p_j\beta_{mj}(ye^{-\bm\gamma^\T {\bm x}})\}^2}
\\&~~~~~~
-\frac{\delta \{\bar{\bm B}_m(y_1e^{-\bm\gamma^\T {\bm x}})-\bar{\bm B}_m(y_2e^{-\bm\gamma^\T {\bm x}})\}^{\otimes 2}}{\{S_m(y_{1}\mid\bm x; \bm\gamma, \bm p)-S_m(y_{2}\mid\bm x; \bm\gamma, \bm p)\}^2},
\end{align}
where $\bm v^{\otimes 2}=\bm v \bm v^\T$  for a column vector $\bm v$, $\bm\beta_m(u)=\{\beta_{m0}(u),\ldots,\beta_{mm}(u)\}^\T$ and  $\bar{\bm B}_m(u)=\{\bar B_{m0}(u),\ldots,\bar B_{mm}(u)\}^\T$.
Denote
\begin{align} \label{eq: deriv of ell(gamma,p) wrt pj for AFT model}
 \Psi_j(\bm\gamma, \bm p) &= \frac{1}{n} \sum_{i=1}^n\frac{\partial \ell_m(\bm\gamma, \bm p;\bm z_i)}{\partial p_j},\quad j=0,\ldots,m.
\end{align}

\begin{theorem}\label{thm: necessary and sufficient condition for AFT model}
For any fixed $\bm\gamma$ suppose $y_{i2}e^{-\bm\gamma^\T {\bm x}_i} \le 1$  for all observed  $(\bm x_i,(y_{i1}, y_{i2}])$ with $y_{i2}<\infty$.
Then  $\tilde{\bm p}(\bm\gamma)$ is a maximizer of $ \ell_m(\bm\gamma, \bm p)$ if and only if
\begin{align}
\label{sufficient-necessary condition 2-interval censoring for AFT model}
\Psi_j \{{\bm\gamma},  \tilde{\bm p}(\bm\gamma)\}&\le 1,
\end{align}
for all $j=0,\ldots,m$ with equality if $\tilde p_j>0$.
\end{theorem}
It is clear that under certain conditions ${\partial^2\ell_m(\bm\gamma,\bm p)}/{\partial \bm p\partial \bm p^\T}$ is a negative and negative definite matrix.
We have fixed-point iteration
\begin{align}
\label{eq: interation for AFT model with general censoring}
    p_j^{[s+1]} &= p_j^{[s]} \Psi_j({\bm\gamma},\bm p^{[s]}),
    \quad\quad\Big. j=0,\ldots,m,\quad s=0,1,2\ldots,\infty.
\end{align}
If $y_{i2}e^{-\bm\gamma^\T {\bm x}_i} \le 1$  for all observed  $(\bm x_i,(y_{i1}, y_{i2}])$ with $y_{i2}<\infty$
then $\Psi_j({\bm\gamma}, {\bm p})\ge 0$ for all $j=0,\ldots,m$
and $\bm p\in \mathbb{S}_{m}$.
Similar to the proof of Theorem 4 of \cite{Peters-and-Walker-1978-siam} we can prove the convergence of ${\bm p}^{[s]}$.
\begin{theorem}\label{thm: convergence of the fixed-point iteration for AFT model}
For any fixed $\bm\gamma$ suppose $y_{i2}e^{-\bm\gamma^\T {\bm x}_i} \le 1$  for all observed  $(\bm x_i,(y_{i1}, y_{i2}])$ with $y_{i2}<\infty$.
If $\bm p^{[0]}$ is in the interior of $\mathbb{S}_{m}$, the sequence $\{\bm p^{[s]}\}$ of (\ref{eq: interation for AFT model with general censoring}) converges to the  maximum approximate profile likelihood estimate (MAPLE) $\tilde{\bm p}(\bm\gamma)$.
\end{theorem}
\paragraph{Algorithm for Finding $(\hat{\bm\gamma},\hat{\bm p})$ for a fixed $m$:}
Let $\tilde{\bm\gamma}$ be an estimate of $\bm\gamma$ such as those proposed by \cite{Jin-etal-2003-bka} and \cite{Tian-and-Cai-bka-2006}.   
\begin{itemize}
\item []
\begin{itemize}
  \item [Step 0:] Start with an initial guess 
$\bm\gamma^{(0)}=\tilde{\bm\gamma}$ of $\bm\gamma$.
    Use (\ref{eq: interation for AFT model with general censoring}) with ${\bm\gamma}= \bm\gamma^{(0)}$,   and the uniform initial $\bm p^{[0]}=\bm u_m\equiv (1,\ldots,1)/(m+1)$  to get $\bm p^{(0)}=\tilde{\bm p}(\bm\gamma)$.
Set $s=0$
  \item [Step 1:] Obtain  $\bm\gamma^{(s+1)}$ with fixed $\bm p=\bm p^{(s)}$ using the Newton-Raphson method starting with $\bm\gamma^{[0]}=\bm\gamma^{(s)}$.
  \item [Step 2:] Choose   ${\bm\gamma}= \bm\gamma^{(s+1)}$.   Then use (\ref{eq: interation for AFT model with general censoring}) with  $\bm p^{[0]}=\bm u_m$ to get $\bm p^{(s+1)}=\tilde{\bm p}(\bm\gamma)$.  Set $s=s+1$.
   \item [Step 3:] Repeat Steps 1 and 2 until convergence. The final $\bm\gamma^{(s)}$ and $\bm p^{(s)}$ are taken as $(\hat{\bm\gamma}, \hat{\bm p})$.
\end{itemize}
\end{itemize}
If $\ell_m(\bm\gamma, \bm p)$ is concave as a function of $\bm\gamma$ then the above algorithm is a point-to-point map and the solution set contains single point.
Convergence  of $(\bm\gamma^{(s)},\bm p^{(s)})$ to $(\hat{\bm\gamma}, \hat{\bm p})$ is guaranteed by the Global Convergence Theorem \citep{Zangwill-1969-book-nonlin-prog}. Proposition \ref{lem: concavity of ell wrt gamma} in Section \ref{sect: asymptotics} suggests that if $n$ is large and $f_m$ is close to $f_0$ then $\ell_m(\bm\gamma, \bm p)$ is  concave with respect to $\bm\gamma$
in a neighborhood of $\bm\gamma_0$.
 \section{Simulation}
We compare the proposed method only with the parametric method for general interval censored data and semparametric competitors whose implementation in R are available such as the rank and the least squares method for right-censored data. 
In all simulation studies, samples of sizes  $n=30, 50, 100$  were generated from Weibull distributions with baseline ($\bm x=\bm 0$) shape 2 and scale 2
according to the AFT model with covariates, $\bm X=(X_1,X_2)$, where $X_1$ and $X_2$ are independent, $X_1$ is uniform$(-1,1)$ and $X_2=\pm 1$ is uniform,  with coefficients ${\bm\gamma}^\T=(\gamma_1,\gamma_2)=(0.5,-0.5)$.
The optimal degrees were chosen from $\{3,\ldots,25\}$ with $\tau_n=12$. Function \verb"ic_par()" of R package \verb"icenReg" \citep{Anderson-Bergman-2017-JSS} was used to obtain parametric maximum likelihood estimates.

In the first simulation study, the proposed method is compared with the parametric method  based on Case $k$  data, where for uncensored censored data $k=0$, for current status data  $k=1$,  and interval censored data with $k$  examinations.   For current status data, the examination time $U$ is uniform$(0,3.66)$ so that $\P(U>T)=50\%$. The general interval censored data with $k$  examinations were generated using the function \verb"simIC_weib()" of  \verb"icenReg" with default arguments when $k=2$, and with \verb"inspections = 5", \verb"inspectLength = 1" when $k=5$.  The censoring probability is 70\% for the  interval censored data with two or more examinations.
 Each sample was used to estimate $\bm\gamma$, $f(\cdot\mid\bm 0)$ and $S(\cdot\mid\bm 0)$ on $[0, \tau_n]$.
In each case, 1000 samples were generated and used to estimate the mean squared errors of the estimates.  The simulation results are shown in Table \ref{tbl: rMSE of gamma-hat and rMISE of f-hat and S-hat in aft model}. We see that,  when data were generated from Weibull distributions, (i)  for small samples, especially the small current status data, the proposed method performs even better than the parametric method in estimating $f(\cdot\mid\bm 0)$ and  $S(\cdot\mid\bm 0)$, (ii)  the parametric method performs a little better than or similar as the proposed method in estimating the regression coefficients for uncensored data or interval censored data with many examinations, (iii) but  for small samples interval censored data with fewer examinations the proposed method performs even better than the parametric method in estimating the regression coefficient.

\begin{table}
  \centering
  \caption{Simulated root mean (integrated) squared errors of estimates of the regression coefficients (the baseline density  and survival functions at $\bm x= \bm 0$) using
  the proposed estimators and the parametric maximum likelihood estimators (in parentheses) based on Case $k$ interval censored data  
  }\label{tbl: rMSE of gamma-hat and rMISE of f-hat and S-hat in aft model}
  \begin{tabular}{crcccc}
& &\multicolumn{4}{c}{Proposed (Parametric)}\\
$k$&$n$ &$\gamma_{1}$&$\gamma_{2}$&$f(\cdot|\bm 0)$&$S(\cdot|\bm 0)$\\
0&$30$ & 0.172 (0.166) & 0.098 (0.096) & 0.114 (0.146) & 0.117 (0.152) \\
 &$50$ & 0.129 (0.128) & 0.077 (0.075) & 0.097 (0.105) & 0.094 (0.113) \\
 &$100$& 0.088 (0.087) & 0.053 (0.052) & 0.079 (0.074) & 0.073 (0.083) \\
1&$30$ & 0.315 (0.428) & 0.199 (0.288) & 0.183 (1.623) & 0.201 (0.319) \\
 &$50$ & 0.232 (0.278) & 0.136 (0.160) & 0.140 (0.324) & 0.141 (0.213) \\
 &$100$& 0.163 (0.191) & 0.097 (0.113) & 0.116 (0.159) & 0.113 (0.137) \\
2&$30$ & 0.225 (0.232) & 0.132 (0.130) & 0.137 (0.195) & 0.145 (0.181) \\
 &$50$ & 0.170 (0.174) & 0.094 (0.094) & 0.120 (0.135) & 0.114 (0.132) \\
 &$100$& 0.111 (0.113) & 0.066 (0.066) & 0.096 (0.091) & 0.085 (0.094) \\
5&$30$ & 0.227 (0.236) & 0.128 (0.128) & 0.156 (0.194) & 0.152 (0.174) \\
 &$50$ & 0.164 (0.169) & 0.101 (0.101) & 0.119 (0.135) & 0.114 (0.134) \\
 &$100$& 0.117 (0.118) & 0.070 (0.068) & 0.097 (0.092) & 0.086 (0.097) \\
  \end{tabular}
\\
The data are uncensored if $k=0$. The censoring rate is 100\% for current status data $(k=1)$, and 70\% for other interval censored data $(k=2,5)$.
\end{table}
In the second simulation, the proposed estimator $\hat{\bm\gamma}$ is compared with  the parametric maximum likelihood estimator
, the rank-based estimator
, and the least squares estimator
, which are implemented in  R package \verb"aftgee" \citep[see][for the details about this package and for more references]{Chiou-et-al-JSS-2014} for the right-censored data with uniform$(0, c)$ right-censoring variable with chosen $c$ to achieve  the specified censoring rates 30\% and 70\%.
From the  results of this simulation given in Table  \ref{tbl: rMSE of estimates of gamma in aft model} we see that the proposed method performs better than the rank and the least squares methods, similar to (even better than) the parametric method if censoring rate is low (high).
\begin{table}
  \centering
  \caption{Simulated root mean squared errors of estimates of the regression coefficients using the proposed, 
  the parametric maximum likelihood
  , the rank-based
  , and the least squares estimators
  }\label{tbl: rMSE of estimates of gamma in aft model}
  \begin{tabular}{crcccccccc}
&   &\multicolumn{2}{c}{Proposed}&\multicolumn{2}{c}{Parametric}
    &\multicolumn{2}{c}{Rank}&\multicolumn{2}{c}{Least squares}  \\
Rate&$n$ &$\gamma_{1}$&$\gamma_{2}$&$\gamma_{1}$&$\gamma_{2}$&$\gamma_{1}$&$\gamma_{2}$
&$\gamma_{1}$&$\gamma_{2}$\\
30\%&$30$ & 0.211 & 0.121 & 0.217 & 0.121 & 0.249 & 0.145 & 0.260 & 0.149\\
    &$50$ & 0.157 & 0.093 & 0.159 & 0.093 & 0.180 & 0.106 & 0.188 & 0.110\\
    &$100$& 0.105 & 0.063 & 0.108 & 0.061 & 0.127 & 0.072 & 0.136 & 0.076\\
70\%&$30$ & 0.319 & 0.187 & 0.404 & 0.612 & 0.482 & 0.525 & 0.469 & 0.528\\
    &$50$ & 0.240 & 0.145 & 0.284 & 0.365 & 0.344 & 0.332 & 0.331 & 0.325\\
    &$100$& 0.162 & 0.105 & 0.178 & 0.136 & 0.217 & 0.166 & 0.216 & 0.153\\
 \end{tabular}
\\
Rate, censoring rate.
\end{table}
\section{Breast Cosmesis Data}
This dataset as described in \cite{Finkelstein-and-Wolfe-1985-biom} and \cite{Finkelstein-1986-biom} is used to study the cosmetic effects of cancer therapy. The time-to-breast-retractions in months ($T$) were  subject to interval censoring and were measured for 94 women among them 46 received radiation only ($X=0$) (25 right-censored, 3 left-censored and 18 interval censored) and 48 received radiation plus chemotherapy ($X=1$)  (13 right-censored, 2 left-censored and 33 interval censored). The right-censored event times were for those women who did not experienced cosmetic deterioration.  The therapy effect on the event time was assessed by many authors. For example, \cite{Hanson-and-Johnson-2004-jcgs-bayesian} fitted the data by a Bayesian AFT model using the mixture of Dirichlet processes \citep{Antoniak-1974-aos-mixture-dirichlet-process} to approximate the baseline survival function and obtained an estimated effect 0.57; \cite{Tian-and-Cai-bka-2006} fitted the data by the AFT model using a Markov chain
Monte Carlo  based resampling method and obtained an estimated effect 0.52 with standard error 0.16; and
\cite{Zhang-and-Davidian-2008-biom-Smooth-SP} used a so-called ``seminonparametric density'' estimator of \cite{Gallant-and-Nychka-1987-econom} and obtained estimated effect $0.95$ with standard error 0.280.

\begin{figure}
\centering
  \makebox{\includegraphics[width=4.7in]{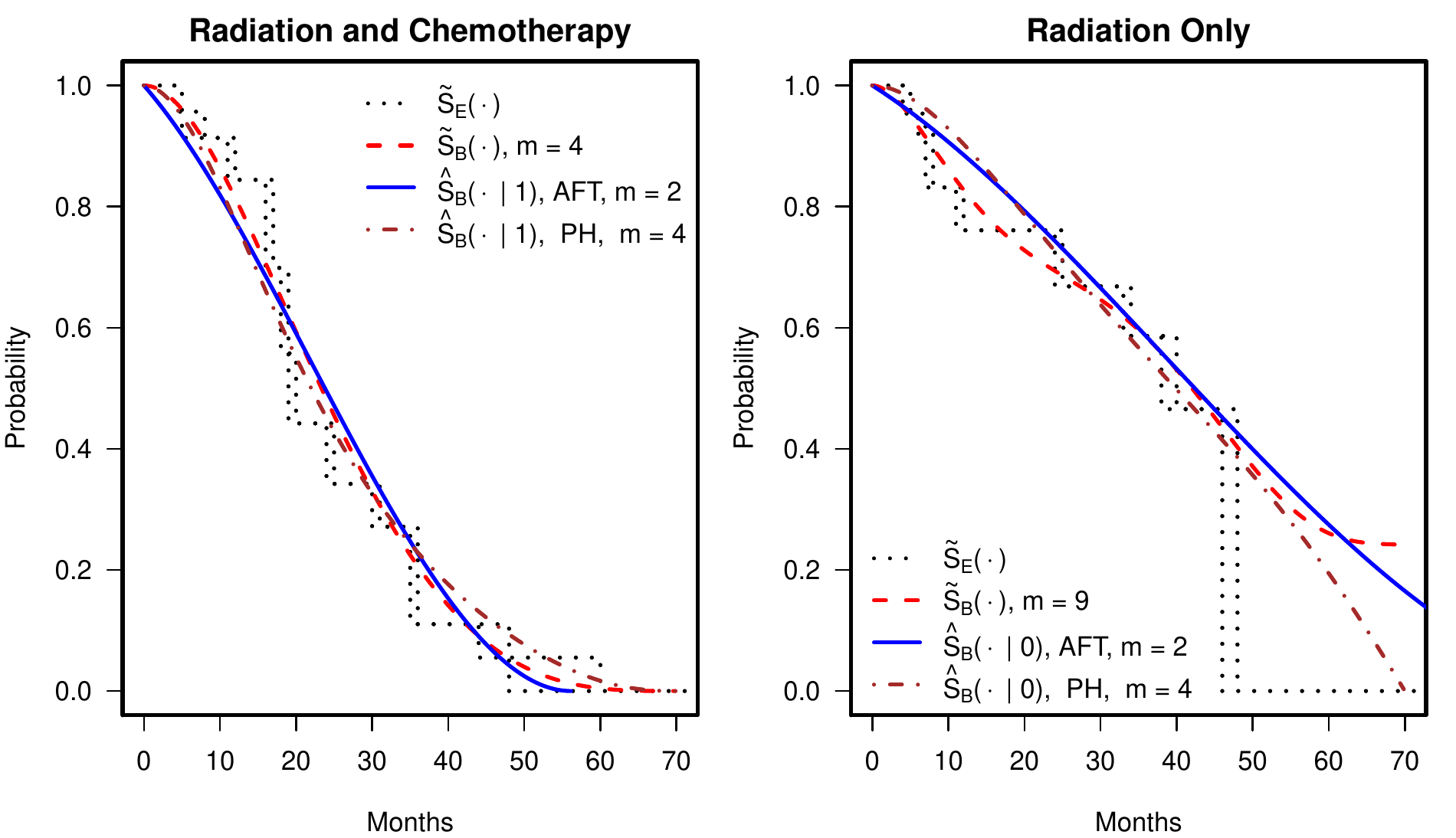}}
\caption{Estimated survival functions for breast cosmetic data:
 the NPMLE $\tilde S_{\mathrm{E}}(\cdot)$ (dotted), and the MABLE $\tilde S_{\mathrm{B}}(\cdot)$ (dashed) from the two samples separately, $\hat S_{\mathrm{B}}(\cdot|\bm x)$  using AFT model (solid) and PH model (dash-dotted) based on combined sample.}\label{fig: Breast-Cosmetic-data}
\end{figure}

The estimated survival curves are shown in Figure \ref{fig: Breast-Cosmetic-data}  where $\tilde S_{\mathrm{E}}$ and $\tilde S_{\mathrm{B}}$ represent, respectively, the NPMLE and the MABLE of $S$ based on each of the two samples, and $\hat S_{\mathrm{B}}(\cdot\mid\bm x)$ represents the proposed estimate based on the combined sample using AFT model  or the estimate of \cite{Guan-arXiv-2019} using PH model.
From this figure it can be seen that, although the two models give similar estimates for radiation and chemotherapy  and for radiation only up to about 45 months, the AFT model fit the data better for months 45 through 60. It is reasonable to believe that the survival probability at month 60 is significantly larger than 0 due to nature and the high percentage of right-censored observations. The AFT model gives an estimated effect 0.572 with standard error 0.123. This is almost the same as the posterior median obtained by \cite{Hanson-and-Johnson-2004-jcgs-bayesian} and close to those in \cite{Tian-and-Cai-bka-2006} but quite different from those given in \cite{Finkelstein-1986-biom}, \cite{Goetghebeur-and-Ryan-2000-biom}, \cite{Betensky-etal-2002-statist-med}, and \cite{Zhang-and-Davidian-2008-biom-Smooth-SP}. The estimated survival curves given by the latter are similar to those of the present paper.
\section{Asymptotic Results}\label{sect: asymptotics}
\setcounter{condition}{-1}
The following assumptions are needed. 
\begin{assump}\label{A1}
The support $\mathcal{X}$ of covariate $\bm X$ is  {compact}  and  $\E({\bm X}{\bm X}^\T)$ is positive definite.
\end{assump}
\begin{assump}\label{A2}
For each  $\tau_n>0$, there exist  $\rho>0$ and ${\bm p}_0=(p_{01},\ldots,p_{0m})^\T$ such that $p_{0i}\ge 0$ ($i=0,\ldots,m$),  $\sum_{i=0}^mp_{0i}=\pi(\bm x_0)$ and,  uniformly in $t\in[0,\tau_n]$,
\begin{equation}\label{approx for f(t|x0)}
\frac{f_m(t; {\bm p}_0)-f_0(t)}{f_0(t)}=\mathcal{O}(m^{-\rho/2}).
\end{equation}
\end{assump}
The positive definiteness of $\E(\bm X\bm X^\T)$ is equivalent to $\Pr(\bm c^\T {\bm X}=0)<1$ for all nonzero $\bm c\in \mathbb{R}^d$. If the right-hand-side of (\ref{approx for f(t|x0)}) is zero for some $m=m_0$, i.e., $f_m(t; {\bm p}_0)=f_0(t)$ for all $t\in[0,\tau_n]$,  then \ref{A2} is true for all $m\ge m_0$ with a zero right-hand-side of (\ref{approx for f(t|x0)}) \citep[see Lemma 2.2 of][]{Guan-2017-jns}.
\begin{prop}\label{lem: concavity of ell wrt gamma}
Suppose that  $\bm\gamma=\bm\gamma_0$,  $f_0$ has continuous second derivative on $[0,\tau_0)$  and   $\E(T^2\mid \bm X=\bm 0)<\infty$ if $\tau_0=\infty$, $\tau_0^2f_0'(\tau_0)+\tau_0f_0(\tau_0)\le 0$ if $\tau_0<\infty$.
 Then, as $n\to\infty$,
$n^{-1}H_n(\bm\gamma,f_0)$ converges almost surely to $H(\bm\gamma,f_0)$ which is negative semi-definite.
Moreover, under Assumption \ref{A1}, if  $\lim_{n\to\infty} n_0/n<1$ or $\lim_{n\to\infty} n_0/n=1$ but $\Pr\{f_0(T)= c_0 /T\mid \bm X=\bm 0\}<1$ for all $c_0>0$, then $H(\bm\gamma,f_0)$ is negative definite.
\end{prop}
\begin{rem}
 Under the conditions of Proposition \ref{lem: concavity of ell wrt gamma}, for $n$ large enough $H_n(\bm\gamma,f_0)$ is almost surely negative definite in a neighborhood of $\bm\gamma_0$.
\end{rem}
\begin{rem}
The condition $\tau_0^2f_0'(\tau_0)+\tau_0f_0(\tau_0)\le 0$ is fulfilled if given $\bm X= \bm 0$ the time $T$ has  a truncated Weibull distribution with shape $\sigma$ and scale $\kappa$ on $[0,\tau_0]$, $W_{\tau_0}(\sigma, \kappa)$, and $\tau_0\ge \kappa$.
\end{rem}

We shall study the large sample property of the proposed estimation under the following conditions regarding to uncensored, Case 1, and Case 2 interval censored data.
\begin{condition}\label{C0}
The event time $T$ is uncensored, $\tau_n\le\tau\le\infty$, $S_0(\tau_n)=\mathcal{O}(n^{-1})$,   and
$|f_0(\theta t)/f_0(t)-1|\le C_0|1-\theta|$ for all $t\in [0,\tau_n]$,  $\theta>0$ and some constant $C_0$ independent of $\bm x_0\in\mathcal{X}$.
\end{condition}
\begin{condition}\label{C1}
The event time $T$ is subject to Case 1 interval censoring, and
given $\bm X=\bm x$ the examination time $U$ has cdf $G(\cdot |  \bm x)$ on $[\tau_l,\tau_u]$, $0<\tau_l<\tau_u<\tau_n\le\tau\le\infty$.
\end{condition}
\begin{condition}\label{C2}
The event time $T$ is subject to Case $2$  interval censoring,
given $\bm X=\bm x$ the {\em observed}  examination times  $\bm U=(U_1, U_2)$ have joint cdf $G(\cdot |  \bm x)$ on $\{(u_1, u_2)\in \mathbb{R}^2: 0<\tau_l\le u_1<u_2\le \tau_u\}$, and $\tau_u<\tau_n\le\tau\le\infty$.
\end{condition}
Under \ref{A1}, Condition
\ref{C0} is satisfied if given $\bm X= \bm 0$ the time $T$ has  a truncated Weibull distribution  $W_{\tau_0}(\sigma, \kappa)$ on $[0,\tau_0]$.
In the following we assume that $\Pr(Y_2 e^{-\bm\gamma^\T\bm X}\le \tau)=1$. We define distance $D_i^2(\bm\gamma,\bm p)$ under Condition $i$, $i=0,1,2$, in the following.
\begin{align}\label{eq: D0^2(gamma, p; x0)}
D_0^2(\bm\gamma,\bm p)&= (\bm\gamma-\bm\gamma_0)^\T\E({\bm X} {\bm X}^\T)(\bm\gamma-\bm\gamma_0)+\chi_0^2(\bm \gamma, \bm p),\\\label{eq: D1^2(gamma, p; x0)}
D_1^2(\bm \gamma, \bm p)
&=
\E\left\{\frac{\{S_m (Ue^{-{\bm\gamma}^\T {\bm X}};\bm p)-S_0(Ue^{-{\bm\gamma}_0^\T {\bm X}})\}^2}{S_0(Ue^{-{\bm\gamma}_0^\T {\bm X}})\{1-S_0(Ue^{-{\bm\gamma}_0^\T {\bm X}})\}}
\right\},\\\label{eq: D2^2(gamma, p; x0)}
D_{2}^2(\bm \gamma, \bm p)&=
\E\left\{\frac{\bm W (\bm U, {\bm X};\bm \theta)^\T
    \bm A(\bm U, {\bm X})
\bm W (\bm U, {\bm X};\bm \theta)}{S_0(U_1e^{-{\bm\gamma}_0^\T {\bm X}})-S_0(U_2e^{-{\bm\gamma}_0^\T {\bm X}})}
\right\},
  \end{align}
where
\begin{align*}
\chi_0^2(\bm \gamma, \bm p)&= \int_{\cal{X}}\int_0^{\tau_0} \left\{\frac{f_m(te^{-(\bm\gamma-\bm\gamma_0)^\T {\bm x}};\bm p)}{f_0(t)}-1\right\}^2 f_0(t) dt dH(\bm x),\\
\bm W(\bm u, {\bm x};\bm \theta)&=\Big\{
        S_m\Big(\frac{u_1}{e^{{\bm\gamma}^\T {\bm x}}};\bm p\Big)-S_0\Big(\frac{u_1}{e^{{\bm\gamma}^\T {\bm x}}}\Big),
        S_m\Big(\frac{u_2}{e^{{\bm\gamma}^\T {\bm x}}};\bm p\Big)-S_0\Big(\frac{u_2}{e^{{\bm\gamma}^\T {\bm x}}}\Big)\Big\}^\T,
\end{align*}
and $\bm A(\bm u, {\bm x})=\{A_{ij}(\bm u, {\bm x})\}$ is a  symmetric matrix with entries
    $A_{11}(\bm u, {\bm x})=
        \{1-S_0(u_2e^{-{\bm\gamma}_0^\T {\bm x}})\}/\{1-S_0(u_1e^{-{\bm\gamma}_0^\T {\bm x}})\}$, $A_{12}(\bm u, {\bm x})=-1$, and $A_{22}(\bm u, {\bm x})={S_0(u_1e^{-{\bm\gamma}_0^\T {\bm x}})}/{S_0(u_2e^{-{\bm\gamma}_0^\T {\bm x}})}$.

We have the following results about the rate of convergence in terms of the above distances.
\begin{theorem}\label{thm: rate of convergence of for aft model}
Let $n_k$ be the number of observations that are subject to Case $k$ censoring  and $\rho_k=\lim_{n\to\infty} n_k/n$, $k=0,1,2$.  Under Assumptions \ref{A1} and \ref{A2},  if Conditions 0, 1, and 2 are fulfilled and $m=C n^{1/\rho}$ for some constant $C$ then, for any $\alpha>1$, $\sum_{k=0}^2 \rho_k D_k^2(\hat{\bm\gamma}, \hat{\bm p})
=\mathcal{O}\{(\log\log n)^\alpha /n\}$, almost surely.  
\end{theorem}
\begin{theorem}\label{thm: induced rates of convergence for aft model}
Under the conditions of Theorem \ref{thm: rate of convergence of for aft model},  if $n_0=n$, then $\chi_0^2(\bm \gamma_0, \hat{\bm p})=\mathcal{O}\{(\log\log n)^\alpha n^{-1+4/\rho}\}$ and $\Vert \hat{\bm\gamma}-\bm\gamma\Vert^2=\mathcal{O}\{(\log\log n)^\alpha /n\}$, almost surely, for any $\alpha>1$;
if $n_0=0$, and
$\tilde{\bm\gamma}$ is an estimate such that $\Vert \tilde{\bm\gamma}-\bm\gamma_0\Vert^2=\mathcal{O}(n^{-1+\epsilon})$ for some $\epsilon>0$, then the MAPLE $\tilde{\bm p}(\tilde{\bm\gamma})$ satisfies
$\sum_{k=1}^2 \rho_k D_k^2(\bm\gamma_0, \tilde{\bm p})=\mathcal{O}(n^{-1+\epsilon'+2/\rho})$ whenever $\epsilon'>\epsilon$.
\end{theorem}
\begin{rem}
  From this theorem it follows that if $\rho$ is large then the convergence rate  of $\sum_{k=1}^2 \rho_k D_k^2(\bm\gamma_0, \tilde{\bm p})$ can be very close to $\mathcal{O}(n^{-1})$ when $\tilde{\bm\gamma}$ is asymptotically normal.
\end{rem}

 \section{Concluding Remarks}\label{sect: concluding remarks}
The proposed approximate likelihood method is even better than some parametric methods based on models with known and fixed dimension due to the lack of robustness of such when sample size is small which is often the case in survival analysis of rare disease and reliability analysis for expensive product. Thus approximated models with unknown dimension enjoy the properties of efficiency and nonparametric robustness.

\def\cprime{$'$} \def\cprime{$'$}
  \def\polhk#1{\setbox0=\hbox{#1}{\ooalign{\hidewidth
  \lower1.5ex\hbox{`}\hidewidth\crcr\unhbox0}}}

\section{Appendix}
\subsection{Proof of Theorem \ref{thm: necessary and sufficient condition for AFT model}}
\begin{proof}
If $y_2e^{-{\bm\gamma}^\T \bm x} \le 1,$ then the negative-definiteness of ${\partial^2\ell_m(\bm\gamma,\bm p;\bm z)}/{\partial \bm p\partial \bm p^\T}$ implies that  $\ell_m(\bm\gamma, \bm p)$ is strictly concave on the compact and convex set $\mathbb{S}_{m}$ for the fixed $\bm\gamma$.
By the optimality condition for convex optimization \citep{Boyd-and-Vandenberghe-book-convex-optimization-2004}  we have that $\tilde{\bm p}$ is the unique maximizer of $\ell_m(\bm\gamma,\bm p)$ if and only if
\begin{equation}\label{sufficient-necessary condition 1-interval censoring}
\nabla_{\bm p}  \ell_m(\bm\gamma,\tilde{\bm p})^\T(\bm p-\tilde{\bm p})\le 0,\quad  \mbox{ for all $\bm p\in \mathbb{S}_{m^*}$},
\end{equation}
where $\nabla_{\bm p}  \ell_m(\bm\gamma, {\bm p})=\partial \ell_m(\bm\gamma,  {\bm p})/\partial \bm p$.
Therefore $\tilde{\bm p}$ is a maximizer of $\ell_m(\bm\gamma, {\bm p})$ for the fixed $\bm\gamma$ if and only if
\begin{align}
\label{sufficient-necessary condition-interval censoring}
\Psi_j({\bm \gamma},\tilde{\bm p})\le & n^{-1}\nabla_{\bm p}  \ell_m(\bm\gamma,\tilde{\bm p})^\T \tilde{\bm p}=1,
\end{align}
for all $j\in\mathbb{I}_0^{m}$ with equality if $\tilde p_j>0$.
The proof is complete.
\end{proof}
\subsection{Proof of Theorem \ref{thm: convergence of the fixed-point iteration for AFT model}}
\begin{proof}
Following the proof of Theorems 1 and 2 and the Corollary of \cite{Peters-and-Walker-1978-siam} we
define
 $\bm \Pi=\diag\{\bm p\}$ and $\bar {\bm \Psi}(\bm p,\bm\gamma)=\bm\Pi {\bm \Psi}(\bm p,\bm\gamma),$
 where ${\bm \Psi}(\bm p,\bm\gamma)=[\Psi_0(\bm p,\bm\gamma),\ldots, \Psi_{m}(\bm p,\bm\gamma)]^\T$. Then
$\bar {\bm \Psi}(\bm p,\bm\gamma)=n^{-1}\bm \Pi \nabla_{\bm p} \ell_m(\bm\gamma,\bm p).$
Its gradient is
$$ \nabla_{\bm p} \bar{\bm \Psi}(\bm p,\bm\gamma)=\frac{\partial \bar{\bm \Psi}(\bm p,\bm\gamma)}{\partial\bm p^\T}
= \frac{1}{n}\diag\{\nabla_{\bm p} \ell_m(\bm\gamma,\bm p)\}+\frac{1}{n}\bm\Pi\frac{\partial \nabla_{\bm p} \ell_m(\bm\gamma,\bm p)}{\partial\bm p^\T}
$$
$$= \frac{1}{n}\diag\Big\{\frac{\partial \ell_m(\bm\gamma,\bm p)}{\partial\bm p}\Big\}+ \frac{1}{n} \bm\Pi\frac{\partial^2 \ell_m(\bm\gamma,\bm p)}{\partial\bm p\partial\bm p^\T}.$$
For any norm on $\mathbb{R}^{m+1}$ we have
$\bar{\bm \Psi}(\bm p,\bm\gamma)-\tilde {\bm p}=\nabla_{\bm p} \bar{\bm \Psi}(\tilde {\bm p},\bm\gamma)(\bm p-\tilde {\bm p})
+\mathcal{O}(\Vert \bm p-\tilde {\bm p}\Vert^2). $
Consider $\nabla_{\bm p} \bar{\bm \Psi}(\tilde {\bm p},\bm\gamma)$ as an operator on subspace
$\mathbb{Z}_m=\{\bm z\in R^{m+1} : \bm 1^\T\bm z=0\}. $
If all components of $\tilde{\bm p}$ are positive then $n^{-1}\nabla_{\bm p} \ell_m(\bm\gamma,\tilde {\bm p})=\bm 1$, and
$\nabla_{\bm p} \bar{\bm \Psi}(\tilde {\bm p},\bm\gamma)=I_{m+1}-\bm Q, $
where $\bm Q=-{n}^{-1}\tilde {\bm \Pi} {\partial^2 \ell_m(\bm\gamma,\tilde {\bm p})}/{\partial\bm p\partial\bm p^\T}.$
From (\ref{eq: 1st deriv of ell_m wrt p for AFT model}) and (\ref{eq: 2nd deriv of ell_m wrt p for AFT model}) it follows that  $\bm Q$ is a left stochastic matrix and $\tilde {\bm p}^\T{\partial^2 \ell_m(\bm\gamma,\tilde {\bm p})}/{\partial\bm p\partial\bm p^\T}=
-{n}^{-1}{\partial \ell_m(\bm\gamma,\tilde {\bm p})}/{\partial\bm p^\T}=- \bm 1^\T$. So
  $\mathbb{Z}_m$ is   invariant under  $\bm Q$.

Define an inner product $\langle\cdot,\cdot\rangle$ by $\langle \bm u,\bm v\rangle=\bm u^\T\tilde {\bm \Pi}^{-1}\bm v$ for $\bm u$, $\bm v$ in $\mathbb{Z}_m$. It can be easily shown that, with respect to this inner product, $\bm Q$ is symmetric and positive   semidefinite  on $\mathbb{Z}_m$:
$$\langle \bm u,\bm Q\bm v\rangle=\bm u^\T\tilde {\bm \Pi}^{-1}\bm Q\bm v=
- \frac{1}{n}\bm u^\T \frac{\partial^2 \ell_m(\bm\gamma,\tilde {\bm p})}{\partial\bm p\partial\bm p^\T}\bm v,
$$
$$\langle \bm Q\bm u,\bm v\rangle=\bm u^\T \bm Q^\T\tilde {\bm \Pi}^{-1}\bm v=
- \frac{1}{n}\bm u^\T \frac{\partial^2 \ell_m(\bm\gamma,\tilde {\bm p})}{\partial\bm p\partial\bm p^\T}\bm v,
$$
$$\langle \bm u,\bm Q\bm u\rangle=
- \frac{1}{n}\bm u^\T \frac{\partial^2 \ell_m(\bm\gamma,\tilde {\bm p})}{\partial\bm p\partial\bm p^\T}\bm u
\ge 0.$$
Let $\lambda_0$ and $\lambda_{m}$ be the smallest and largest eigenvalues of $\bm Q$ associated with eigenvectors in $\mathbb{Z}_m$. Then the operator norm of $\nabla_{\bm p} \bar{\bm \Psi}(\tilde {\bm p},\bm\gamma)$ on $\mathbb{Z}_m$ w.r.t. this inner product equals
$\max\{|1-\lambda_0|,|1-\lambda_{m}|\}$. It is clear that $0\le \lambda_0\le\lambda_{m}\le 1$ because $\bm Q$ is a left stochastic matrix. Because ${\partial^2 \ell_m(\bm\gamma,\tilde {\bm p})}/{\partial\bm p\partial\bm p^\T}<0$  we have $\lambda_0>0$.
Similar to the proof of Theorem 2 of \cite{Peters-and-Walker-1978-siam} the assertion of theorem follows. If $\tilde{\bm p}$ contains zero component(s),
say $\tilde p_j=0$, $j\in J_0$, deleting the $j$-th row  and $j$-th column  of the vectors and matrices in the above proof for all $j\in J_0$ we can show that  the iterates $p_j^{(s)}$, $s\in\ints(0,\infty)$, converge to $\tilde p_j$ as $s\to \infty$ for all $j\notin J_0$. Because $\sum_{j=0}^{m}  p_j^{(s)}=1$ and $ p_j^{(s)}\ge 0$, $j\in\ints(0,m)$, for those $j\in J_0$,  $p_j^{(s)}$ converges to zero as $s\to \infty$.
The proof of Theorem \ref{thm: convergence of the fixed-point iteration for AFT model} is complete.
\end{proof}
\subsection{Proof of Proposition \ref{lem: concavity of ell wrt gamma}}
\begin{proof}
The derivatives of $\ell(\bm\gamma,f_0;\bm z)$ w.r.t. $\bm\gamma$ are
\begin{align}\nonumber
 \frac{\partial\ell(\bm\gamma,f_0;\bm z)}{\partial \bm \gamma}&=-(1-\delta)\left\{ 1 +\frac{y e^{-\bm\gamma^\T {\bm x}} f_0'(ye^{-\bm\gamma^\T {\bm x}} )}{f_0( ye^{-\bm\gamma^\T {\bm x}})}\right\} {\bm x}
 \\\label{eq: 1st derivative of ell wrt gamma} &~~~
+\delta e^{-\bm\gamma^\T {\bm x}}\frac{y_1f_0(y_{1}e^{-\bm\gamma^\T {\bm x}})-y_2f_0(y_{2}e^{-\bm\gamma^\T {\bm x}})}{S_0(y_{1}e^{-\bm\gamma^\T {\bm x}})-S_0(y_{2}e^{-\bm\gamma^\T {\bm x}})} {\bm x},
\\\non
  \frac{\partial^2\ell(\bm\gamma,f_0;\bm z)}{\partial \bm\gamma\partial \bm\gamma^\T}&=
  (1-\delta)\left\{\frac{y e^{-\bm\gamma^\T {\bm x}} f_0'(ye^{-\bm\gamma^\T {\bm x}} )}{f_0( ye^{-\bm\gamma^\T {\bm x}})}\right.
   \\\nonumber &~~~
   \left.+ \frac{y^2 e^{-2\bm\gamma^\T {\bm x}} [ f_0(ye^{-\bm\gamma^\T {\bm x}} )f_0''( ye^{-\bm\gamma^\T {\bm x}})-\{f_0'( ye^{-\bm\gamma^\T {\bm x}})\}^2]}{f_0^2( ye^{-\bm\gamma^\T {\bm x}})}\right\} {\bm x}{\bm x}^\T
    \\\nonumber
 &~~~
-\delta e^{-\bm\gamma^\T {\bm x}}\left[\frac{y_1f_0(y_{1}e^{-\bm\gamma^\T {\bm x}})-y_2f_0(y_{2}e^{-\bm\gamma^\T {\bm x}})}{S_0(y_{1}e^{-\bm\gamma^\T {\bm x}})-S_0(y_{2}e^{-\bm\gamma^\T {\bm x}})} \right.\\\non
&~~~
+ e^{-\bm\gamma^\T {\bm x}}\left.\frac{\{y_1f_0(y_{1}e^{-\bm\gamma^\T {\bm x}})-y_2f_0(y_{2}e^{-\bm\gamma^\T {\bm x}})\}^2}{\{S_0(y_{1}e^{-\bm\gamma^\T {\bm x}})-S_0(y_{2}e^{-\bm\gamma^\T {\bm x}})\}^2}\right.\\\label{eq: 2nd derivative of ell wrt gamma}
&~~~
\left.+e^{-\bm\gamma^\T {\bm x}}\frac{y_1^2f_0'(y_{1}e^{-\bm\gamma^\T {\bm x}})-y_2^2f_0'(y_{2}e^{-\bm\gamma^\T {\bm x}})}{S_0(y_{1}e^{-\bm\gamma^\T {\bm x}})-S_0(y_{2}e^{-\bm\gamma^\T {\bm x}})}\right] {\bm x} {\bm x}^\T.
\end{align}
If $\tau_0=\infty$, then $\E(T\mid \bm X=\bm x)
=\int_0^\infty t e^{-\bm\gamma^\T {\bm x}} f_0(te^{-\bm\gamma^\T {\bm x}})dt<\infty$ implies $\lim_{y\to\infty}yf_0(y)=0$.
Similarly, $\E(T^2\mid \bm X=\bm x)<\infty$ implies  $\lim_{y\to\infty}y^2f_0(y)=0$ and  thus
 $\int_0^\infty t^2 e^{-2\bm\gamma^\T {\bm x}} f'_0(te^{-\bm\gamma^\T {\bm x}})dt =-2\E(T\mid \bm X=\bm x).$
The latter ensures   $\lim_{y\to\infty}y^2f_0'(y)=0$.
Therefore we have
\begin{equation}\label{two limits}
\lim_{y\to\infty}y^{k+1}f_0^{(k)}(y)=0, \quad k=0,1.
\end{equation}
By (\ref{eq: 2nd derivative of ell wrt gamma}) and the SLLN we have, as $n\to\infty$,
$$  \frac{1}{n}\frac{\partial^2\ell(\bm\gamma,f_0)}{\partial \bm\gamma\partial \bm\gamma^\T} \to H(\bm\gamma,f_0)=
   -\sum_{k=0}^2 \rho_k \int_{\mathcal{X}} I_k(\bm x) e^{-2\bm\gamma^\T {\bm x}}{\bm x}{\bm x}^\T dH(\bm x),\quad a.s.,$$
where $\rho_k$ is the probability that $T$ is Case $k$ censored so that $\rho_0+\rho_1+\rho_2=1$, and
 $I_k(\bm x)$ ($k=0,1,2$) are given below.

For Case 0 data, using substitution $t=ye^{-\bm\gamma^\T {\bm x}}$ we have
\begin{align*}
  I_0(\bm x)
&= e^{2\bm\gamma^\T {\bm x}} \int_0^{\tau_0} \Big[t\frac{\{f_0'(t)\}^2}{f_0(t)}- f_0'(t)-tf_0''(t)\Big]  tdt \\
& =e^{2\bm\gamma^\T {\bm x}}  \int_0^{\tau_0} \Big\{ \frac{t f_0'(t)}{f_0(t)}+1\Big\}^2 f_0(t)dt-e^{2\bm\gamma^\T {\bm x}} \{\tau_0^2f_0'(\tau_0)+\tau_0f_0(\tau_0)\}.
    \end{align*}
It follows from (\ref{two limits}) if $\tau_0=\infty$ and the assumption $\tau_0^2f_0'(\tau_0)+\tau_0f_0(\tau_0)\le 0$ if $\tau_0<\infty$  that
$$I_0(\bm x)\ge e^{2\bm\gamma^\T {\bm x}}  \int_0^b \Big\{ \frac{t f_0'(t)}{f_0(t)}+1\Big\}^2 f_0(t)dt$$ for all $\bm x\in\mathcal{X}$.
Therefore $\Pr\{f_0(T)= c_0 /T\mid \bm X=\bm 0\}<1$, for all $c_0>0$, implies that  $I_0(\bm x)>0$ for all $\bm x\in\mathcal{X}$.

For Case 1 data, 
\begin{align*}
I_1(\bm x)&=   \E \left(e^{\bm\gamma^\T {\bm X}}\left[\frac{{-Uf_0(Ue^{-\bm\gamma^\T {\bm X}})}}{1-S_0(Ue^{-\bm\gamma^\T {\bm X}})} + e^{-\bm\gamma^\T {\bm X}}\frac{\{Uf_0(Ue^{-\bm\gamma^\T {\bm X}})\}^2}{\{1-S_0(Ue^{-\bm\gamma^\T {\bm X}})\}^2}  \right.\right.\\
&\quad\quad
\left.\left. +e^{-\bm\gamma^\T {\bm X}}\frac{{ -U^2f_0'(Ue^{-\bm\gamma^\T {\bm X}})}}{1-S_0(Ue^{-\bm\gamma^\T {\bm X}})}\right]  {I(0<T<U) } \,\Big|\,\bm X=\bm x\right)\\
&~~~ +   \E \left(e^{\bm\gamma^\T {\bm X}}\left[\frac{{Uf_0(Ue^{-\bm\gamma^\T {\bm X}})}}{S_0(Ue^{-\bm\gamma^\T {\bm X}})} + e^{-\bm\gamma^\T {\bm X}}\frac{\{Uf_0(Ue^{-\bm\gamma^\T {\bm X}})\}^2}{\{S_0(Ue^{-\bm\gamma^\T {\bm X}})\}^2}  \right.\right.\\
&\quad\quad
\left.\left. +e^{-\bm\gamma^\T {\bm X}}\frac{{U^2f_0'(Ue^{-\bm\gamma^\T {\bm X}})}}{S_0(Ue^{-\bm\gamma^\T {\bm X}})}\right] {I(U<T<\infty) } \Big|\bm X=\bm x\right)\\
&=   \E \left( \left[ \frac{\{Uf_0(Ue^{-\bm\gamma^\T {\bm X}})\}^2}{ 1-S_0(Ue^{-\bm\gamma^\T {\bm X}}) }
  +    \frac{\{Uf_0(Ue^{-\bm\gamma^\T {\bm X}})\}^2}{S_0(Ue^{-\bm\gamma^\T {\bm X}})}   \right]\Big|\bm X=\bm x  \right)\\
&=   \int_{0}^{\infty}  \frac{ \{uf_0(ue^{-\bm\gamma^\T {\bm x}})\}^2}{ \{1-S_0(ue^{-\bm\gamma^\T {\bm x}})\} S_0(ue^{-\bm\gamma^\T {\bm x}})} dG(u\mid\bm x).
\end{align*}
Similarly, for Case $2$  data, we have 
\begin{align*}
I_2(\bm x)
 &= \mathop{\iint}_{u_1<u_2}  \left[ \frac{\{u_1f_0(u_1e^{-\bm\gamma^\T {\bm x}})\}^2}{ 1-S_0(u_1e^{-\bm\gamma^\T {\bm x}}) }+\frac{\{u_1f_0(u_1e^{-\bm\gamma^\T {\bm x}})-u_2f_0(u_2e^{-\bm\gamma^\T {\bm x}})\}^2}{S_0(u_1e^{-\bm\gamma^\T {\bm x}})-S_0(u_2e^{-\bm\gamma^\T {\bm x}}) }
 \right.\\
  &\quad\quad\quad\left.+    \frac{\{u_2f_0(u_2e^{-\bm\gamma^\T {\bm x}})\}^2}{S_0(u_2e^{-\bm\gamma^\T {\bm x}})}   \right]dG(u_1,u_2\mid \bm x)\\
  &= \mathop{\iint}_{u_1<u_2} \left[ \frac{\bm \psi(\bm u, {\bm x}; \bm\gamma)^\T\bm A(\bm u, {\bm x})\bm \psi(\bm u, {\bm x}; \bm\gamma)}{S_0(u_1e^{-\bm\gamma^\T {\bm x}})-S_0(u_2e^{-\bm\gamma^\T {\bm x}}) }
\right]dG(\bm u\mid \bm x),
    \end{align*}
where
$\bm \psi(\bm u, {\bm x}; \bm\gamma)=\{u_1f_0(u_1e^{-\bm\gamma^\T {\bm x}}), u_2f_0(u_2e^{-\bm\gamma^\T {\bm x}})\}^\T$ and $\bm A(\bm u, {\bm x})$ is the same as in (\ref{eq: D2^2(gamma, p; x0)}).
Clearly $I_k(\bm x)>0$ for $k=1,2$ and all $\bm x\in\mathcal{X}$.
The negative-definiteness of $H(\bm\gamma,f_0)$ follows.
   The proof is complete.
\end{proof}
\subsection{Proof of Theorem \ref{thm: rate of convergence of for aft model}}
\begin{proof} Assume that the data are arranged so that $\delta_i=k$ for $i=n_{k-1}+1,\ldots,n_k$, $k=0,1,2$, where $n_{-1}=0$.
  Let $\bm\gamma_0$ be the true value of $\bm\gamma$. Define
   $\tilde{\mathcal{R}}({\bm\gamma},\bm p)=\{\ell({\bm\gamma}_0,f_0)-\ell_m({\bm\gamma},\bm p)\}/n$. By (\ref{eq: approx bernstein loglikelihood for AFT model}) and Taylor expansion we obtain
  \begin{align}\label{eq: R(gamma,p)}
    \tilde{\mathcal{R}}({\bm\gamma},\bm p)
&=-\sum_{0\le i\le j\le 1}\tilde{\mathcal{R}}_{ij}({\bm\gamma},\bm p)
+\frac{1}{2}\sum_{i=0}^1\tilde{\mathcal{R}}_{i2}({\bm\gamma},\bm p)
+ \sum_{i=0}^1o\{\tilde{\mathcal{R}}_{i2}({\bm\gamma},\bm p)\},
  \end{align}
where $\tilde{\mathcal{R}}_{00}({\bm\gamma},\bm p) =-n^{-1}\sum_{i=0}^{n_0} (\bm\gamma-\bm\gamma_0)^\T {\bm x}_i$,
  \begin{align*}
\tilde{\mathcal{R}}_{0k}({\bm\gamma},\bm p)&= \frac{1}{n}\sum_{i=0}^{n_0} \Biggr\{\frac{f_m(y_ie^{-\bm\gamma^\T {\bm x}_i};\bm p)}{f_0(y_ie^{-\bm\gamma_0^\T {\bm x}_i})}-1\Biggr\}^k, \quad k=1,2,\\
\tilde{\mathcal{R}}_{1k}({\bm\gamma},\bm p)&=\frac{1}{n}\sum_{i=n_0+1}^n \Biggr\{ \frac{S_m(y_{1i}e^{-\bm\gamma^\T {\bm x}_i};\bm p)-S_m(y_{2i}e^{-\bm\gamma^\T {\bm x}_i};\bm p)}{S_0(y_{1i}e^{-\bm\gamma_0^\T {\bm x}_i}) -S_0(y_{2i}e^{-\bm\gamma_0^\T {\bm x}_i})}-1\Biggr\}^k, \quad k=1,2.
  \end{align*}
By Stirling formula, for all real $x$,
\begin{equation}\label{eq: expansion of exp(x)-1}
\left| e^{x}
-\sum_{i=0}^j\frac{1}{i!}
x^i\right|\le \frac{e^2\sqrt{j+1} |x|^{j+1}e^{|x|}}{\sqrt{2\pi}(j+1)!},\quad j=0,1,2,\ldots.
\end{equation}
All the large sample statements in the following proofs are {\em almost sure}.

If $\rho_0>0$, under Condition \ref{C0}, by the LIL and Kolmogorov's SLLN 
we have
  \begin{align}\non\label{eq: R00}
    \tilde{\mathcal{R}}_{00}({\bm\gamma},\bm p)
&=-\rho_0(\bm\gamma-\bm\gamma_0)^\T \E({\bm X})
+\mathcal{O}[\{(\bm\gamma-\bm\gamma_0)^\T\Var({\bm X})(\bm\gamma-\bm\gamma_0)\}^{1/2}\\
&~~~~~~\cdot (\log\log n/n)^{1/2}],
\\\non
    \tilde{\mathcal{R}}_{02}({\bm\gamma},\bm p)
     &= \rho_0\int_{\cal{X}}\int_0^{\tau_0 e^{\bm\gamma_0^\T\bm x}} e^{-{\bm\gamma}_0^\T {\bm x}} \Biggr\{\frac{f_m(ye^{-\bm\gamma^\T {\bm x}};\bm p)}{f_0(ye^{-\bm\gamma_0^\T {\bm x}})}-1\Biggr\}^2
f_0(ye^{-\bm\gamma_0^\T {\bm x}})dy dH(\bm x)\\\non
&~~~+o\left\{\chi_0^2(\bm \gamma, \bm p)\right\}\\\label{eq: R02}
&=\rho_0\chi_0^2(\bm \gamma, \bm p)+o\left\{\chi_0^2(\bm \gamma, \bm p)\right\},
  \end{align}
  where
  \begin{align*}
\chi_0^2(\bm \gamma, \bm p)&= \int_{\cal{X}}\int_0^{\tau_0 e^{\bm\gamma_0^\T\bm x}} e^{-{\bm\gamma}_0^\T {\bm x}} \frac{\{f_m(ye^{-\bm\gamma^\T {\bm x}};\bm p)- f_0(ye^{-\bm\gamma_0^\T {\bm x}})\}^2}{f_0(ye^{-\bm\gamma_0^\T {\bm x}})} dy dH(\bm x)\\
&= \int_{\cal{X}}\int_0^{\tau_0 e^{\bm\gamma_0^\T\bm x_0}} \left\{\frac{f_m(te^{-(\bm\gamma-\bm\gamma_0)^\T {\bm x}};\bm p)}{f_0(t)}-1\right\}^2 f_0(t) dt dH(\bm x).
  \end{align*}
By the LIL we have
  \begin{align}\non
    \tilde{\mathcal{R}}_{01}({\bm\gamma},\bm p)
     &=\rho_0\int_{\cal{X}}\int_0^{\tau_0 e^{\bm\gamma_0^\T\bm x}} e^{-{\bm\gamma}_0^\T {\bm x}}\Biggr\{ \frac{f_m(ye^{-\bm\gamma^\T {\bm x}};\bm p)}{f_0(ye^{-\bm\gamma_0^\T {\bm x}})}-1 \Biggr\}
f_0(ye^{-\bm\gamma_0^\T {\bm x}})dy dH(\bm x)\\\non
&~~~+\mathcal{O}\left\{\chi_0(\bm \gamma, \bm p) (\log\log n/n)^{1/2}\right\} \\\label{eq: R01}
     &=\rho_0\E\{e^{(\bm\gamma-\bm\gamma_0)^\T {\bm X}}-1\}  +\mathcal{O}\{\chi_0(\bm \gamma, \bm p) (\log\log n/n)^{1/2}\}.
  \end{align}
  By (\ref{eq: expansion of exp(x)-1}) we get
\begin{align}\non
\tilde{\mathcal{R}}_0({\bm\gamma},\bm p)&\equiv-\sum_{i=0}^1\tilde{\mathcal{R}}_{0i}({\bm\gamma},\bm p)+\frac{1}{2}\tilde{\mathcal{R}}_{02}({\bm\gamma},\bm p)+
o\{\tilde{\mathcal{R}}_{02}({\bm\gamma},\bm p)\}\\\non
&=\E\{e^{(\bm\gamma-\bm\gamma_0)^\T {\bm X}}-1-(\bm\gamma-\bm\gamma_0)^\T {\bm X}\}
 +\frac{1}{2}\chi_0^2(\bm \gamma, \bm p)\\\non
&~~~~ +\mathcal{O}\{\chi_0(\bm \gamma, \bm p) (\log\log n/n)^{1/2}\}\\\non
&~~~~+\mathcal{O}[\{(\bm\gamma-\bm\gamma_0)^\T\Var({\bm X})(\bm\gamma-\bm\gamma_0)\}^{1/2} (\log\log n/n)^{1/2}]\\\non
&~~~~~~+o\left\{\chi_0^2(\bm \gamma, \bm p)\right\}\\\non
&=\frac{1}{2}D_0^2(\bm\gamma,\bm p)+\mathcal{O}\{\chi_0(\bm \gamma, \bm p) (\log\log n/n)^{1/2}\}\\\non
&~~~~+\mathcal{O}[\{(\bm\gamma-\bm\gamma_0)^\T\Var({\bm X})(\bm\gamma-\bm\gamma_0)\}^{1/2} (\log\log n/n)^{1/2}]\\\label{expansion of R(gamma,p) under condition 0}
&~~~~+\frac{1}{6}\mathcal{O}[\E\{ e^{|(\bm\gamma-\bm\gamma_0)^\T {\bm X}|}|(\bm\gamma-\bm\gamma_0)^\T {\bm X}|^3\}]+o\left\{\chi_0^2(\bm \gamma, \bm p)\right\}.
\end{align}
Under Conditions \ref{C1} and \ref{C2}, if $\rho_k>0$, $k=1,2$, similar to the proof of Proposition \ref{lem: concavity of ell wrt gamma}, 
we have $ \tilde{\mathcal{R}}_{11}({\bm\gamma},\bm p)
 = \mathcal{O}\{\sum_{k=1}^2\rho_k D_{k}(\bm \gamma, \bm p)(\log\log n/n)^{1/2}\}$ by LIL and 
$\tilde{\mathcal{R}}_{12}({\bm\gamma},\bm p)= \sum_{k=1}^2\rho_k D_{k}^2(\bm \gamma, \bm p)+o\{\sum_{k=1}^2\rho_k D_{k}^2(\bm \gamma, \bm p)\}$ by Kolmogorov's SLLN. Thus
\begin{align}\non
\tilde{\mathcal{R}}_1({\bm\gamma},\bm p)&\equiv- \tilde{\mathcal{R}}_{11}({\bm\gamma},\bm p)+\frac{1}{2}\tilde{\mathcal{R}}_{12}({\bm\gamma},\bm p)+
o\{\tilde{\mathcal{R}}_{12}({\bm\gamma},\bm p)\}\\\non
&=\frac{1}{2} \sum_{k=1}^2\rho_k D_{k}^2(\bm \gamma, \bm p)
 +\mathcal{O}\left\{\sum_{k=1}^2\rho_k D_{k}(\bm \gamma, \bm p)(\log\log n/n)^{1/2}\right\}\\\label{expansion of R(gamma,p) under conditions 1 and 2}
  &\quad +o\left\{\sum_{k=1}^2\rho_k D_{k}^2(\bm \gamma, \bm p)\right\}.
\end{align}
If $\sum_{k=0}^2\rho_k D_{k}^2(\bm \gamma, \bm p)= (\log\log n)^{\alpha}/n$ then
$$\tilde{\mathcal{R}}({\bm\gamma},\bm p)=\frac{1}{2}(\log\log n)^{\alpha}/n
+\mathcal{O}\{(\log\log n)^{(\alpha+1)/2}/n\}+o\{(\log\log n)^{\alpha}/n\}.$$
While at $(\bm\gamma,\bm p)=(\bm\gamma_0,\bm p_0)$, if $m=C n^{1/\rho}$
then, by  \ref{A2},  $\sum_{k=0}^2\rho_k D_{k}^2(\bm \gamma, \bm p)=\mathcal{O}(1/n)$ and
$\tilde{\mathcal{R}}({\bm\gamma},\bm p)=\mathcal{O}(1/n) +o(1/n)$.
Thus  $D_0^2(\hat{\bm \gamma}, \hat{\bm p})\le (\log\log n)^{\alpha}/n$  if $\alpha>1$.
\end{proof}
\subsection{Proof of Theorem \ref{thm: induced rates of convergence for aft model}}
\begin{proof} Theorem \ref{thm: rate of convergence of for aft model} implies that  $\Vert \hat{\bm\gamma}-\bm\gamma_0\Vert^2 =\mathcal{O}\{(\log\log n)^{\alpha}/n\}$ under Assumption \ref{A1} and  Condition 0. By the reverse triangular inequality we have
   \begin{align}\label{2nd estimate of chi0(gamma,p) under condition 0}
\chi_0^2(\bm \gamma, \bm p)
& \ge \{\chi_0(\bm \gamma_0, \bm p)-\tilde\chi_0(\bm \gamma, \bm p)\}^2.
  \end{align}
  where $$\tilde\chi_0^2(\bm \gamma, \bm p)=\int_{\cal{X}}\int_0^{\tau_0 e^{\bm\gamma_0^\T\bm x}} e^{-{\bm\gamma}_0^\T {\bm x}} \frac{\{f_m(ye^{-\bm\gamma^\T {\bm x}};\bm p)- f_m(ye^{-\bm\gamma_0^\T {\bm x}};\bm p)\}^2}{f_0(ye^{-\bm\gamma_0^\T {\bm x}})} dy dH(\bm x).$$
Noting that $|f_m''(t;\bm p)|\le C_2 m^2$, Theorem \ref{thm: rate of convergence of for aft model} with $i=0$, together with
(\ref{2nd estimate of chi0(gamma,p) under condition 0}), imply that $\chi_0^2(\bm \gamma_0,\hat{\bm p})=\mathcal{O}\{(\log\log n)^{\alpha}n^{-1+4/\rho}\}$.
%

Under Condition 1,
by trianglar inequality
\begin{align} \label{est of D1^2}
  D_1^2(\bm \gamma, \bm p) & \le 2D_1^2(\bm \gamma_0, \bm p) + 2\tilde D_1^2(\bm \gamma, \bm p),
\end{align}
where
$$\tilde D_1^2(\bm \gamma, \bm p)=\int_{\cal{X}} \int_{\tau_l}^{\tau_u}
\frac{\{S_m(ue^{-{\bm\gamma}^\T {\bm x}};\bm p)-S_m(ue^{-{\bm\gamma}_0^\T {\bm x}};\bm p)\}^2}{S_0(ue^{-{\bm\gamma}_0^\T {\bm x}})\{1-S_0(ue^{-{\bm\gamma}_0^\T {\bm x}})\}} dG(u\mid\bm x)  dH(\bm x).$$
We also have
\begin{align} \label{lower bound of D1^2}
  D_1^2(\bm \gamma, \bm p) & \ge \{D_1(\bm \gamma_0, \bm p)
- \tilde D_1(\bm \gamma, \bm p)\}^2.
\end{align}
Noting that $f_m(t;\bm p)\le (m+1)/\tau_n$ we have,
\begin{align}\non \label{upper bound of tilde D1^2}
\tilde D_1^2(\bm \gamma, \bm p)
&\le  \frac{(m+1)^2}{\tau_n^{2}}
\int_{\cal{X}} \int_{\tau_l}^{\tau_u}
\frac{u^2 (e^{-{\bm\gamma}_0^\T{\bm x}}-e^{-{\bm\gamma}^\T{\bm x}})^2dG(u\mid\bm x)}{S_0(ue^{-{\bm\gamma}_0^\T{\bm x}})\{1-S_0(ue^{-{\bm\gamma}_0^\T{\bm x}})\}}   dH(\bm x)\\
 &= \tau_n^{-2} (m+1)^2 \mathcal{O}\Big\{(\bm\gamma-\bm\gamma_0)^\T \bm M_0(\bm\gamma-\bm\gamma_0)\Big\},
\end{align}
where
$$\bm M_0=\int_{\cal{X}} \int_{\tau_l}^{\tau_u}
\frac{u^2 e^{-2{\bm\gamma}_0^\T {\bm x}}{\bm x} {\bm x}^\T dG(u\mid\bm x) }{S_0(ue^{-{\bm\gamma}_0^\T {\bm x}})\{1-S_0(ue^{-{\bm\gamma}_0^\T {\bm x}})\}}   dH(\bm x).$$
Because $\mathcal{X}$ is compact and $0<\tau_l<\tau_u<\infty$,
$$\mathcal{O}\{(\bm\gamma-\bm\gamma_0)^\T \bm M_0(\bm\gamma-\bm\gamma_0)\}=\mathcal{O}(\Vert \bm\gamma-\bm\gamma_0\Vert^2).$$
For any fixed $\bm\gamma\in \mathbb{B}_{n,\epsilon}$,
if $m = C n^{1/\rho}$ and  $D_1^2(\bm \gamma_0, \bm p)= n^{-1+\epsilon'+2/\rho}$ for any $\epsilon'>\epsilon$,
we can show that there exists an $\eta>0$ so that
$\tilde{\mathcal{R}}({\bm\gamma},\bm p)\ge \eta n^{-1+\epsilon'+2/\rho}$, a.s..  While at $\bm p=\bm p_0$ we have, by (\ref{est of D1^2}), $\tilde{\mathcal{R}}({\bm\gamma},\bm p)= \mathcal{O}(n^{-1+\epsilon+2/\rho})$, a.s.. Thus for any fixed $\bm\gamma\in \mathbb{B}_{n,\epsilon}$, the maximizer $\bm p(\bm\gamma)$ of $\ell_m({\bm\gamma},\bm p)$ satisfies $D_1^2\{\bm \gamma_0, \bm p(\bm\gamma)\} \le n^{-1+\epsilon'+2/\rho}$ whenever $\epsilon'>\epsilon$.

Under Condition 2, 
$D_{2}^2(\bm \gamma, \bm p)  \le 2\tilde D_{2}^2(\bm \gamma, \bm p)+2D_{2}^2(\bm \gamma_0, \bm p)$ and $D_2^2(\bm \gamma, \bm p)  \ge \{D_2(\bm \gamma_0, \bm p)
- \tilde D_2(\bm \gamma, \bm p)\}^2$,
where
\begin{align*}
\tilde D_{2}^2(\bm \gamma, \bm p)&=  \E\left\{\frac{\tilde{\bm W } (\bm U, {\bm X};\bm \theta)^\T
    \bm A(\bm U, {\bm X})
\tilde{\bm W } (\bm U, {\bm X};\bm \theta)}{S_0(U_1e^{-{\bm\gamma}_0^\T {\bm X}})-S_0(U_2e^{-{\bm\gamma}_0^\T {\bm X}})}
\right\},\\
\tilde{\bm W }(\bm u, {\bm x};\bm \theta)&=\!\Big\{
        S_m\Big(\frac{u_1}{e^{{\bm\gamma}^\T {\bm x}}};\bm p\Big)\!-\!S_m\Big(\frac{u_1}{e^{{\bm\gamma}_0^\T {\bm x}}};\bm p\Big),
        S_m\Big(\frac{u_2}{e^{{\bm\gamma}^\T {\bm x}}};\bm p\Big)\!-\!S_m\Big(\frac{u_2}{e^{{\bm\gamma}_0^\T {\bm x}}};\bm p\Big)\Big\}^\T.
\end{align*}
The proof under Condition 2 is similar to that under Condition 1.
\end{proof}

 \end{document}